\documentclass[12pt]{amsart}
\usepackage{amsmath, amssymb, amsfonts}
\makeatletter
\def\ps@pprintTitle{%
   \let\@oddhead\@empty
   \let\@evenhead\@empty
   \let\@oddfoot\@empty
   \let\@evenfoot\@oddfoot
}
\makeatother

\newtheorem {Theorem}{Theorem}[section]

\newtheorem {Proposition}[Theorem]{Proposition}

\newtheorem {theorem}{Theorem}[section]
\newtheorem {lemma}[theorem]{Lemma}

\newtheorem {remark}[theorem]{Remark}

\newcommand{\di}{\mbox{div}}

\newcommand{\R}{\mathbb{R}}
\newcommand{\p}{\partial}
\newcommand{\ba}{\begin{array} }
\newcommand{\ea}{\end{array} }

\def\R{{\mathbb{R}}}

\newcommand{\be}{\begin }

\numberwithin{equation}{section}

\def\ba{\begin{aligned}}
\def\ea{\end{aligned}}
\def\be{\begin{equation}}
\def\ee{\end{equation}}
\def\leqs{\lesssim}
\def\d#1{\dot{#1}}
\def\da{\f{D}{Dt}}
\def\t#1{\tilde{#1}}
\def\eps{\varepsilon}
\def\ts{\tau}

\def\r{\rho}

\def\div{\text{div}}

\def\O{\Omega}

\def\f{\frac}
\def\p{\partial}
\def\|{\parallel}
\def\a{\alpha}
\def\g{\nabla}
\def\ld{\lambda}

\def\u{\mathbf{u}}
\def\v{\delta}

\def\O{\Omega}

\allowdisplaybreaks

\begin{document}
	\vspace{1cm}
\title[Global large strong solutions to 2d CNS]{Global large strong solutions to the radially symmetric compressible Navier-Stokes equations in 2D solid balls}
\author{Xiangdi Huang, Mengluan SU, Wei YAN$^*$, Rongfeng Yu}
\thanks{* Corresponding author.}

\address{Xiangdi HUANG\hfill\break\indent
	Institute of Mathematics,
\hfill\break\indent
	Academy of Mathematics and Systems Sciences,\hfill\break\indent
	Chinese Academy of Sciences, Beijing 100190, China}
\email{xdhuang@amss.ac.cn}
\address{Mengluan SU\hfill\break\indent
	Institute of Applied Physics and Computational Mathematics, Beijing, 100088, P.R.China\break\indent
	and \hfill\break\indent
	School of Mathematics and Statistics,\hfill\break\indent
	Northeast Normal University, Changchun, 130022, P.R.China}
\address{Wei YAN\hfill\break\indent
	School of Mathematics,\hfill\break\indent
	Jilin Unversity, Changchun 130012, China}
\email{wyanmath@jlu.edu.cn}
\address{Rongfeng YU\hfill\break\indent
	School of Mathematics,\hfill\break\indent
	Jilin Unversity, Changchun 130012, China}
\email{rfyu@jlu.edu.cn}
\maketitle

\begin{abstract}

In this paper, we consider the initial-boundary value problems  of the compressible isentropic Navier-Stokes equations with density-dependent viscosity on two dimensional solid balls which was first introduced by Kazhikhov where shear viscosity $\mu$ is assumed to be constant and the bulk viscosity $\lambda$ is a polynomial of density up to power $\beta$.

Under the condition of $\beta>1$, we prove the global existence of the radially symmetric strong solutions to the Kazhikhov models under Dirichlet boundary conditions for arbitrary large initial smooth data. Moreover, the density is shown to be uniformly bounded with respect to time when $\beta\in (\max\{1,\frac{\gamma+2}{4}\},\gamma]$. This improves the previous result of \cite{2016Huang,2022Huang} for general 2D domains where they require $\beta>4/3$ to ensure global existence and is the first result concerning the global existence of classical solutions to the radially symmetric compressible Navier-Stokes equations in 2D solid balls under Dirichlet boundary condition.

\end{abstract}

\vspace*{10pt}

\section{Introduction and main results}

This paper is concerned with the compressible isentropic Navier-Stokes system in a bounded domain $\O\subset \R^2$: 
\begin{equation} \label{eq1}
\left\{\begin{array}{lr} 
\rho_t +\di(\rho \u)=0,\\
(\rho \u)_t+\di(\rho \u \otimes \u)+ \nabla P=\mu\Delta\u+\nabla ((\mu+\lambda)\di \u),
\end{array}\right.
\end{equation}
with the initial data
\begin{equation}\label{initial data}
\rho(x,0)=\rho_{0}(x), \u(x,0)=\u_{0}(x),~~~~x\in \O,
\end{equation}
and Dirichlet boundary condition
\begin{equation}\label{BD}
\u(x,t)=0~~~~~on ~~\p\O.
\end{equation}
Here $x\in\O$ is the spatial coordinate, $t\geq 0$ is time. $\rho \geq 0$, $\u$ and $P$ denote the fluid density, velocity and pressure, respectively.

Whether the 2D compressible isentropic Navier-Stokes equations with constant viscosity will admit global classical solutions with arbitrary large initial smooth data remains a longstanding open problem. The major breakthrough is due to Merle(\cite{merle1, merle2}) where they proved that the equations will develop finite time singularity for multi-dimensional radially symmetric solutions when the domain is a solid ball. Thus the only expected and avaliable result is to prove the result for annulus which is already well known in(\cite{Choe2005}).

It turns out that the equations become more interesting and flexible when the viscosity is assumed to be density-dependent, especially when Vaigant and Kazhikhov introduced their model with $\mu=const.$, $\lambda=\rho ^\beta $ in periodic domains. Indeed, under the condition $\beta>3$, Kazhikhov etcs proved the following surprising reslut: the large strong solutions is shown to be exist globally with arbitrary large initial smooth data. Their work revealed a deep structure balancing the triple terms: the convection term, the diffusion term and the viscosity. Thus it is very interesting and important to investigate their model. This is the main purpose of this paper.

Without lossing of generality, we make the following assumptions.

$\mu$ and $\lambda$ are the viscosity coefficients satisfing the following restrictions:
\begin{equation}
\mu= {\mbox const> 0},\; \lambda(\rho)=b\rho^{\beta},
\end{equation}
with $b>0, ~\beta>1.$  We consider the polytropic gases for which the equation of state is given by 
\begin{equation*}\label{pressure}
P(\rho)=A\rho^{\gamma},
\end{equation*}
where $A$ is a positive constant, $\gamma>1$ is the adiabatic exponent. 

A great deal of mathematicians have made great efforts and achieved fruitful results for the local and global existence of solutions to compressible Navier-Stokes equations. The one-dimensional problem has been studied extensively, see \cite{2011Ding,Liliang2016,li2008vanishing,Qinyao2010,yang2002compressible,Yangyao2001} and references cited therein. In multi-dimensional case, Matsumura and Nishida proved the global existence of classical solutions for compressible Navier-Stokes equations in \cite{matsumura1983initial-boundary}, where they required that the initial data close to a equilibrium state far away from vacuum. When considering general data, one has to face the possible appearance of vacuum. However, as observed in \cite{xin1998blowup,xin2013on}, the smooth solutions to the full compressible Navier-Stokes equations with constant viscosities will blow up in finite time. Some other related results can be found
in \cite{Ch,rozanova2008blow,bian2019finite,duan2019finite-time,2019liwang,yu} and the references therein. When it comes to barotropic compressible Navier-Stokes equations with constant viscosity coefficients, Lions  made a breakthrough in \cite{Lion}, where he proved the global existence of weak solutions for any initial data containing vacuum, as soon as the initial energy was finite and $\gamma\geq \gamma_n (\gamma_n=\dfrac{3n}{n+2}, n=2,3)$. Jiang and Zhang extended Lion's result to $\gamma>1$ for spherically symmetric initial data in \cite{jiang2001on}.  Feireisl, Novotn\'y and Petzeltova \cite{feireisl2001on} improved the condition to $\gamma >3/2$ in three-dimensional space. Huang, Li and Xin showed the global existence of classical solutions to three-dimensional barotropic compressible Navier-Stokes equations for initial data with small total energy but possible large oscillations and containing vacuum states in \cite{huang2012global}. Later, Li and Xin \cite{li2019global} also proved the global existence in two-dimensional space.

When the viscosity coefficients are both dependent of density, Bresch and Desjardins \cite{BD1} proposed a new entropy inequality (BD-entropy) under an additional constraint on the viscosity coefficients, which played an important role in proving the existence of weak solutions. Based on this conclusion, there are some results about weak solution for compressible Navier-Stokes equations. Bresch, Desjardins and G\'{e}rard-Varet showed the existence of global weak solutions for compressible Navier-Stokes equations with an additional appropriate constraints on the initial density profile and domain curvature in \cite{bresch2007on}. Guo, Jiu and Xin proved the global existence of weak solutions for the spherically symmetric initial data in \cite{guo2008spherically}. Li and Xin proved the global existence of weak solutions for two-dimensional and three-dimensional Cauchy problem of barotropic compressible Navier-Stokes equations in \cite{jing2015global}. In the same time, Vasseur and Yu gave the proof of the global existence of weak solutions for three-dimensional Navier-Stokes equations by using a different method in \cite{vasseur2016existence}. A nature question is: Can we improve the regularity of the weak solutions above? Li, Pan and Zhu investigated the local existence of regular solutions for compressible barotropic Navier-Stokes equations with density-dependent viscosities in \cite{Li2016Recent,2015On}. Luo and Zhou extended the result in \cite{2019luozhou}. Recently, Xin and Zhu \cite{xin2019global} proved the global-in-time well-posedness of regular solutions for a class of smooth initial data for Cauchy problem. So it is important to investigate the classical solutions for multi-dimensional compressible Navier-Stokes equations with degenerate viscosities whether exist globally.

For the case that the viscosity coefficients depend on the density and degenerate at the vacuum, there are more interesting phenomena. Such model was derived from the Boltzmann equations by Liu, Xin and Yang in \cite{liu1998}. 

Another interesting model was introduced by Vaigant and Kazhikhov in 1995 where they\cite{1995Kazhikhov} proved that there exists a unique classical solution for two-dimensional barotropic compressible Navier-Stokes equations with $\mu=const.$, $\lambda=\rho ^\beta (\beta>3)$ in periodic domains when the initial density is away from vacuum. In 2016, by employing a new structure and regularity criterion, also combined with Zlotnik's inequality, the condition on $\beta$ was weakened to $\beta>4/3$ , Huang and Li  \cite{2016Huang} also established the global existence of classical solutions to the periodic domain in presence of vacuum. Recently, by carelly choosen weight on time and space of the initial data, this was generalized to 
the whole space \cite{2022Huang} with vacuum at far field. See also Jiu-Wang-Xin \cite{JWX} for the Cauchy problem with positive density at far field.

It is worth noting that the above methods in \cite{2016Huang,2022Huang} depends largely on some commutator's estimates which works well for both periodic domain or the whole space but failed to bounded domains with Dirichlet boundary conditions. Thus it becomes very interesting and chanllenging to derive similar results for bounded domains. This is the main goal we want to address in this paper. The simpliest case is no doubt radially symmetric domains.

In this paper, we consider the global existence and large time behavior of spherically symmetric strong solutions to the initial boundary value problem for the compressible Navier-Stokes system \eqref{eq1}. Let $\O$ be a ball of radius $R$ centered at the origin in $\R^2$, and $A$ and $b$  be unity for convenience. We are concerned with spherically symmetric strong solutions to system $\eqref{eq1}$ with initial data $\eqref{initial data}$ and Dirichlet boundary condition $\eqref{BD}$. To this end, we denote 
\begin{equation}\label{sym}
	|x|=r,~~\rho(x,t)=\rho(r,t),~~\u(x,t)=u(r,t)\dfrac{x}{r}.
\end{equation}
The equations  $\eqref{eq1}-\eqref{BD}$ is then transformed to
\begin{equation} \label{eq2}
	\left\{\begin{array}{lr} 
		\rho_t +(\rho u)_r+\dfrac{\rho u}{r}=0,\\
		(\rho u)_t+(\rho u^2)_r+\dfrac{\rho u^2}{r}+ ({\rho^\gamma})_r=[(2\mu+\rho^\beta)(u_r+\dfrac{u}{r})]_r,
	\end{array}\right.
\end{equation}
with the initial data
\begin{equation}\label{initial data2}
	\rho(r,0)=\rho_{0}(r), u(r,0)=u_{0}(r),~~~~0\leq r<R,
\end{equation}
and Dirichlet boundary condition
\begin{equation}\label{BD2}
	u(0,t)=u(R,t)=0,~~~~~t>0.
\end{equation}

The main results are described as follows:
\begin{theorem}\label{them1}
	Suppose that 
	\be\label{cdt}
	 \beta>1,~ \gamma>1,
	 \ee 
	 and the initial data $(\rho_0,\u_0)$ satisfy 
	\be \label{HD0}
	\rho_0\in W^{1,q}(\O),\quad \u_0\in H^1(\O),\quad q>2,
	\ee
	and satisfy the following compatibility condition
	~\be\label{HD}
	-\mu\Delta \u_0 -\nabla ((\mu+\lambda(\rho_0)) \ div \u_0) + \nabla P(\rho_0) = \sqrt{\rho_{0}} g,
	\ee
	where ~$g\in L^2(\Omega).$ Then the problem \eqref{eq1}--\eqref{BD} has a unique strong solution satisfying 
	\be
	\rho\in C([0,T); W^{1,q}),\quad \u\in C([0,T);  H^1),
	\ee  and ~\be\left\{
	\ba
	&\|\sqrt{\rho} \u\|_{L^\infty L^2} + \|\rho\|_{L^\infty L^\gamma} + \|\g \u\|_{L^2 L^2} \leq C,\\
	&\|\sqrt{\rho}\dot \u\|_{L^\infty L^2} + \|\g\dot \u\|_{L^2L^2}\leq C.
	\ea\right.\ee
\end{theorem}

A few remarks are in order:
\begin{remark}
In their paper \cite{2016Huang}, one of the the main ingredients is making use of some commutator's estimates combined with the density equation and the variable coefficients for either periodic domains or whole space but failed to apply to bounded domains. It becomes main obstacles to deal with trouble terms arsing from the corresponding commutator, i.e, some boundary terms in both space and time. In our paper, under some new observation and analysis, we successfully bound the boundary term by carefully choosen multipliers on the momentum equation. See Lemma \ref{FR} and below. 
\end{remark}

\begin{remark}
It is worth nothing that Li-Zhang\cite{Lhlz2016} proved an analogous result of Theorem \ref{them1} for free boundary conditions. They can benefit a lot from the commutator's estimates under such a boundary condition. To be more specific,  the extra trouble term $F(R,t)$ in \eqref{frrr} disappears under free boundary condition and thus it is much easier to derive the upper bound of the density. How to treat $F(R,t)$ is the main  obtacle in our analyis which requires much effort and new observations on the structure of the equations.
\end{remark}

\begin{remark}
In Propositions 3.1-3.2, by establishing a higher intergrabity of the density and velocity, we can prove the upper bound of the density depending on the time which gurantees the strong solution can be extended to a global one.
\end{remark} 

Moreover, the density is shown to be uniformly bounded with respect to time when $\beta\in(\max\{1,\frac{\gamma+2}{4}\},\gamma]$. This is also completely new compared to \cite{Lhlz2016,2016Huang,2022Huang}.

\begin{theorem}\label{them2}
Let 
\be
\gamma>1,~\max\{1,\frac{\gamma+2}{4}\}<\beta\le\gamma,
\ee
then there is a positive constant C depending only on $\mu,\beta,\gamma,\|\rho_0\|_{L^\infty}$ and $\|\u_0\|_{H^1}$ such that
\be
\sup_{0\le t<\infty}\|\rho\|_{L^\infty}\le C
\ee
and the following asympotic behaviors holds
\be
\lim_{t\rightarrow\infty}(\|\rho-\rho_s\|_{L^p}+\|\nabla \u\|_{L^p})=0
\ee
for any $p\in[1,\infty)$, where
\be
\rho_s=\frac{1}{|\O|}\int_{\Omega}\rho_0dx.
\ee
\end{theorem}
A few remarks are in order:
\begin{remark}
It seems no hope to establish time-independent upper bound of the density following the proof of Theorem \ref{them1}. Thus we employ Zlotnik's inequality which works well in \cite{2016Huang} etc. By a careful decomposition of the combined equation, it is very lucky to find that the corresponding term contributes exactly in suitable position at a price that $\beta\in (\max\{1,\frac{\gamma+2}{4}\},\gamma]$. This new restriction characterizes a different phynomenon when compared to the Cauchy problem where $\beta>\max\{1,\frac{\gamma+2}{4}\}$. The main reason is due to the control of boundary term involved in the commutator's estimates in deriving the upper bound of the density.
\end{remark} 

\begin{remark}
In the case $\gamma\in (1,2]$, the restiction on $\beta$ is reduced to $1<\beta\le\gamma$.
\end{remark}

\section{Preliminaries}
When the initial density is strictly away from vacuum, the well-known local existence theory can be found in \cite{SALVI, SOL}, which could be stated as follows:
\begin{lemma}\label{local}
	Assume that  $(\rho_0,\u_0)$ satisfy 
	\[
	\inf\limits_{x\in \O}\rho_0(x) >0,\ \rho_0\in H^2,\u_0\in H^2,\u\Big|_{\p\O} = 0.
	\] 
	Then there exists a small time $T>0$ and a constant $C_0>0$ depending only on 
	$\O, \gamma, \beta, \mu, \ld, \|(\rho_0, \u_0)\|_{H^2},\inf\limits_{x\in \O}\rho_0(x)$ such that there exists a unique strong solution $(\rho, \u)$ to the problem \eqref{eq1}--\eqref{BD} in $\O\times(0, T)$ satisfying 
	\[
	\inf\limits_{(x,t)\in \O\times(0,T)}\rho(x,t) > C_0.
	\]
\end{lemma}

\begin{lemma}\cite{Kazhikhov}\label{ME}For any~$ q >2$, there is a constant $C_q>0,$ such that for any function $\u(x)$, if~$\u\big|_{\p\O}=0$ or~$\int_\O \u dx = 0$, then
	\be
	\| \u\|_{q} \leq C_q \parallel \nabla \u \parallel _{\frac{2q}{q+2}}.
	\ee
\end{lemma}

\begin{lemma}\label{lem2}
	Assume that~$\u(x,t)$ and~$u(r,t)$ satisfy \eqref{sym} and \eqref{BD2}, then it holds that
	~\be
	\parallel u \parallel_{\infty} \leq \parallel \div \u \parallel_{2}\leq\parallel \nabla \u \parallel _{2}.
	\ee
\end{lemma}
\proof
For any $r\in (0,R),$ directly calculations lead to
\[\ba
|u(r)|^2 &= 2\int_0^r u\p_r u dr \leq 2\int_0^R|\p_r u||\f{u}{r}| r dr \\
&\leq \int_0^R|\p_r u|^2 r dr + \int_0^R|\f{u}{r}|^2 r dr = \int_0^R \Big(\p_r u + \f{u}{r}\Big)^2 rdr
=\|\div \u\|_2^2,
\ea\]
which implies that~$ \parallel u \parallel _{\infty} \leq \parallel \div \u \parallel_{2}.$
\endproof

The following Sobolev inequality will be used frequently.
\begin{lem}\label{lem3} (See \cite{AD, Kazhikhov}) If~$\u|_{\p\O}=0$ or~$\int_\O \u dx = 0$, then for any ~$q>2$, we have
	~\be\label{Sob1}
	\|\u\|_q\leq C\|\u\|_2^{\f{2}{q}} \|\g \u\|_2^{1-\f{2}{q}},\quad \|\u\|_q \leq C(\|\u\|_2 + \|\g \u\|_2),
	\ee
	and
	~\be\label{Sob}
	\|\u\|_\infty\leq C\|\u\|_2^{\f{q-2}{2q-2}} \|\g \u\|_q^{\f{q}{2q-2}}.
	\ee
\end{lem}

The following Beale--Kato--Majda type inequality can be found in \cite{BKM, hlx}.
\begin{lem}(See \cite{BKM, hlx})
	For any~$2< q < \infty$, there is a constant~$C(q)$, such that the following estimate holds for all ~$\nabla \u \in W^{1,q}(\Omega)$,
	\be\label{CL}
	\parallel \nabla \u \parallel_{\infty} \leq C(\parallel \div \u \parallel_{\infty}+\parallel curl \u \parallel_{\infty} )\log(e + \parallel \nabla ^{2} \u \parallel_{L^q}) + C\parallel \nabla \u \parallel_{L^2} + C.
	\ee
\end{lem}
    
\begin{lem} For any smooth function $f$, we have
	\be\label{CJ}
	\int_\O\rho \f{D}{Dt} f dx = \f{d}{dt}\int_\O \rho f dx,
	\ee
where $\rho$ is the density satisfying $\eqref{eq1}_1$ and  $\frac{D}{Dt}f$ denotes the material derivative of a function f, i.e.
\be
\frac{D}{Dt}f=f_t+\u\cdot\nabla f.
\ee
\end{lem}

The following Zlotnik inequality will be used to get the time-independent upper bound of the
density.
\begin{lem}\label{zlemma} (\cite{zlot}) Let the function $y\in W^{1,1}(0, T)$ satisfys
\[
y'(t) = g(y) + h'(t) \text{ on } [0, T],\quad y(0) = y_0,
\]
with $g\in C(R)$ and $h\in W^{1,1}(0, T)$. If $g(\infty)=-\infty$ and
\be
h(t_2) - h(t_1) \leq N_0 + N_1(t_2 - t_1)
\ee
for all $0 \le t_1 < t_2 \le T$ with some $N_0 \ge 0$ and $N_1 \ge 0$, then
\[
y(t) \le \max(y_0, \zeta_0) + N_0 < \infty \text{ on } [0, T],
\]
where $\zeta$ is a constant such that $g(\zeta) \leq -N_1$ for $\zeta\ge\zeta_0$.
\end{lem}

\section{A priori estimates (I): upper bound of the density}

Let $T>0$ be a fixed time and $(\rho,\u)$ with $\u(x,t)=u(r,t)\dfrac{x}{r}$ be a strong solution to the problem \eqref{eq1}--\eqref{BD} on $\O\times(0, T)$. Then we will give some necessary a priori estimates for $(\rho,\u)$. We begin with the following basic energy estimates.
\begin{lem}\label{EE} {\bf (Basic energy estimates)}
	There is a constant $C>0$, depending only on the initial data $(\rho_{0},  \u_{0})$, such that
	\begin{equation}\label{EE2}
		\sup_{t \in [0,\ T]} (\parallel \sqrt{\rho }\u\parallel_{2}^{2} + \parallel \rho \parallel_{\gamma} ^{\gamma}) + \int _{0}^{T}(\mu\parallel \nabla \u\parallel_{2}^{2} + \parallel(\mu+\lambda(\rho))^{\frac{1}{2}} \ div \u\parallel_{2}^{2} )  dt \leq C.
	\end{equation}
\end{lem}
\proof
It is easy to get the result by taking the inner product of $\eqref{eq1}_2$ with $\u$. We omit the details here.
\endproof

Define
~\be
\theta (\rho)= 2\mu \ln \rho+\frac{1}{\beta} (\rho^{\beta}  -1),
\ee
and the effective viscous flux
~\be\label{EVis}
F = (2\mu+\ld)\div \u - P.
\ee
By making use of the original mass equation $\eqref{eq2}_1$, we have
~\be\label{EDW}
\theta_{t} + u\partial_{r} \theta + F + P = 0.
\ee
Next, thanks to the definition of the effective viscous flux~\eqref{EVis}, we thus rewrite the momentum equation~$\eqref{eq2}_2$ as:
\be\label{XM}
(\rho u )_{t} + \partial_{r} (\rho u^2) + \frac{ 1}{r} \rho u^2 = \partial_{r}F.
\ee
Integrating over~$(R, r)$ gives
$$\frac{d}{dt} \int_{R}^{r} \rho u ds + \rho u^2 +\int_{R}^{r}\frac{\rho u^2}{s} ds = F(r,t)- F(R,\ t).$$
Set
~\be\label{xieta}
\xi=\int_R^r\rho u dr,\quad \eta = \rho u^2 + \int_R^r\f{\rho u^2}{s}ds.
\ee
One has
$$\xi_{t} + \eta - F + F(R,\ t)=0.$$
Adding the above equation to \eqref{EDW} gives
~\be\label{thxi}
(\theta +\xi )_{t} + u \cdot \partial_{r}(\theta+\xi) +\int_{R}^{r} \frac{\rho u^2}{s} ds +P + F(R,\ t)=0.
\ee

Indeed, \eqref{thxi} is exacting the correponding form as (3.42) in \cite{2016Huang} in radially symmetric case which is the key commutator's estimates combined with the density equation and variable coefficients.

Note that $F(R,t) \neq 0$ due to the Dirichlet boundary condition \eqref{BD}. This is one of the main obstacles different from either periodic or the whole space or the paper \cite{Lhlz2016}. In the following, we should give some estimates on the term $F(R,t)$:
\begin{lemma} \label{FR}
It holds that
	\be\label{frrr}
	F(R,\ t) = \frac{1}{R^2} \Big[ \frac{d}{dt} \int_{0}^{R} \rho u r^2 dr + 2\int_{0}^{R} F(r, t) r dr - \int_{0}^{R} \rho u^2 r ds\Big],
	\ee
	and there is a constant $C_1 >1$,such that for all $0\leq t\leq T$,
	\be
	\Big|\int_{0}^{t} F(R, s) ds \Big|\leq C_1 + C\int_0^t\|\rho\|_{\beta}^{\beta}d\ts.
	\ee
\end{lemma}

\proof Multiplying~\eqref{XM} by~$r^2$ gives
$$(r^2 \rho u)_{t}+\partial_{r}(r^2 \rho u^2) - r \rho u ^2 = r^2\partial _{r}F,$$
and then integrating over~$(0, R)$ leads to
$$\int_{0}^{R}(r^2\rho u)_{t} dr + \int_{0}^{R} \partial_{r}(r^2\rho u^2) dr -\int_{0}^{R} r \rho u^2 dr=\int_{0}^{R} r^2\partial_{r}Fdr.$$
Observing that
$$\int_{0}^{R} \partial_{r}(r^2\rho u^2) dr =r^2\rho u^2 \big|_{0}^{R} = 0,$$
and 
$$\int_{0}^{R} r^2\partial_{r}Fdr = r^2 F \big |_{0}^{R}-\int_{0}^{R} F\cdot 2rdr =R^2F(R,\ t)- 2\int_{0}^{R} F\cdot rdr,$$
one has
$$F(R,\ t)= \frac{1}{R^2}\Big [ \frac{d}{dt}\int_{0}^{R} r^2 \rho u dr+ 2\int_{0}^{R} F\cdot rdr-\int_{0}^{R} r \rho u^2 dr\Big].$$

Next, we come to estimate $F(R,t)$.\\
Integating $F(R,t)$ over~$(0,\ t)$, we have
\begin{eqnarray*}
	\int_{0}^{t} F(R, \ts) d\ts &=&  \frac{1}{R^2} \Big[\int_{0}^{R} \rho u r^2 dr \big|_{\ts=t} - \int_{0}^{R} \rho u r^2 dr \big|_{\ts=0} + 2\int_{0}^{t}\int_{0}^{R} F(r, \ts) r dr d\ts\\
	& - &\int_{0}^{t}\int_{0}^{R} \rho u^2 r dr d\ts \Big],
\end{eqnarray*}
and then
\begin{eqnarray*}\label{iFR}
	\begin{split}
		\Big|\int_{0}^{t} F(R,\ \ts) d\ts \Big| & \leq  \frac{1}{R^2} \Big[ R \cdot \int_{0}^{R} \rho |u| r dr + R \int_{0}^{R} \rho_{0}|u_{0}| r dr + 2\int_{0}^{t}\int_{0}^{R} F(r, \ts) r dr d\ts \\
		& \quad+  \int_{0}^{t}\int_{0}^{R} \rho u^2 r dr d\ts \Big]\\
		& \leq  C(1+ T + \int_{0}^{T}\int_{0}^{R} |F(r, \ts)| r dr d\ts  ),
	\end{split}
\end{eqnarray*}
where we used the fact that~$|u| \leq 1+|u|^2$, and $\rho |u| r \leq \rho r + \rho |u|^2 r$.\\
It remains to estimate the term~$\int_{0}^{T}\int_{0}^{R} |F(r, \ts)| r dr d\ts $.
Notice that $$F = (2\mu + \lambda) \div u -P,$$ so we have
\begin{eqnarray*}
	2\pi\int_{0}^{R} |F(r, \ts)| r dr &=& \int_{\Omega} |F(r, \ts)| dx  =  \int_{\Omega}|(2\mu + \lambda) \div u -P| dx\\
	& \leq & \int_{\Omega}(2\mu + \lambda)|\div u| dx +\int_{\Omega}P dx\\
	&  \leq &\int_{\Omega}(2\mu + \lambda)|\div u| dx + C\\
	& \leq & C + \int_{\Omega} (2\mu + \lambda)|\div u|^2 dx + \|\rho\|_{\beta}^{\beta}.
\end{eqnarray*}

Hence we obtain that
\begin{eqnarray*}
	\Big|\int_{0}^{t} F(R,\ \ts) d\ts \Big| & \leq & C(1+ T + \int_{0}^{T}\int_{0}^{R} |F(r, \ts)| r dr d\ts  )\\
	& \leq & C(1+ \int_{0}^{T} \int_{\Omega} (2\mu + \lambda)|\div u|^{2} dx d\ts +C\int_0^t\|\rho\|_{\beta}^{\beta}d\ts)\\
	& \leq & C_1 + C\int_0^t\|\rho\|_{\beta}^{\beta}d\ts,
\end{eqnarray*}
which completes the proof.
\endproof
\begin{remark}\label{R32} 
	We could not get the estimate $\|F(R,t)\|_{L^1(0,T)}$ directly since   $F(R,t)$ contains the term
	\[
	\frac{d}{dt} \int_{0}^{R} \rho u r^2 dr.
	\]
\end{remark}	


Next, we will give the High-order integrability of the density and velcity.\\
Set
~\be\label{Ge}
G(t)\triangleq\int_{0}^{t} F(R, s)ds + 2C_1 + C\int_0^T\|\rho\|_{\beta}^{\beta}d\ts,
\ee
one gets
$$1 \leq C_1 \leq G(t) \leq 3 C_1+ C\int_0^T\|\rho\|_{\beta}^{\beta}d\ts.$$
Then \eqref{thxi} could be rewritten as
~\be\label{thxiG}
(\theta + \xi + G(t))_{t} + u\cdot \nabla(\theta+\xi +G(t)) +\int_{R}^{r} \frac{\rho u^2}{s} ds +P(r, t) =0.
\ee

\begin{Proposition}\label{BJ}
	If~$\beta >1,\gamma>1$, then there exists a constant $T>0$ denpending only on ~$\beta, \gamma, T$ and the initial data, such that
	\be
	\sup_{t\in[0,\ T]}\int_{\Omega} \rho^{2\beta\gamma+1 }dx \leq C.
	\ee
\end{Proposition}

\proof
Let~$f=(\theta + \xi + G(t))_{+}$,
multiplying the \eqref{thxiG} by ~$\rho f^{2\gamma -1}$,
and integrating over $\Omega$, we have
\begin{eqnarray*}
	\int_\O \rho f^{2\gamma -1}\cdot ((\theta &+& \xi + G(t))_{t} + u\cdot \nabla(\theta+\xi +G(t)))dx\\
	&=&  \frac{1}{2\gamma} [\int_\O\rho \cdot\partial_{t}(f^{2\gamma} + \rho u\cdot \partial_{r}(f^{2\gamma})dx]\\
	&=& \frac{1}{2\gamma} [\int_\O\partial_{t}(\rho f^{2\gamma}) - \rho_{t}\cdot f^{2\gamma} + \div (\rho u\cdot f^{2\gamma}) - f^{2\gamma} \cdot \div(\rho u)]\\
	&=&\frac{1}{2\gamma}[\int_\O \partial_{t}(\rho f^{2\gamma}) + \div (\rho u\cdot f^{2\gamma})dx].
\end{eqnarray*}
Notice that $$\int_{\O} \div (\rho u\cdot f^{2\gamma})dx = \int_{\partial \Omega} \rho u\cdot f^{2\gamma}\cdot n ds = 0.$$
We have
\begin{eqnarray*}
	\frac{1}{2\gamma} \frac{d}{dt}\int_{0}^{R} &\rho& f^{2\gamma} dx \leq \int_{\Omega} (\int_{r}^{R} \frac{\rho u^2}{s} ds) (\rho f^{2\gamma-1})dx.\\
	& \leq & C\parallel \rho^{\frac{1}{2\gamma}} f \parallel_{L^{2\gamma}(\Omega)}^{2\gamma-1} \parallel \rho \parallel_{L^{2\beta\gamma+1}(\Omega)}^{\frac{1}{2 \gamma}} \parallel\int_{r}^{R} \rho u^2 \frac{1}{s} ds\parallel_{L^{\frac{2\beta\gamma+1}{\beta}}(\Omega)}\\
	&\leq & C\parallel \rho^{\frac{1}{2\gamma}} f \parallel_{L^{2\gamma}(\Omega)}^{2\gamma-1} \parallel \rho \parallel_{L^{2\beta\gamma+1}(\Omega)}^{\frac{1}{2 \gamma}} \parallel \rho u^2 \frac{x}{r^2} \parallel_{L^{\frac{2(2\beta \gamma +1)}{2\beta\gamma+1+2\beta}}(\Omega)}\\
	&\leq & C\parallel \rho^{\frac{1}{2\gamma}} f \parallel_{L^{2\gamma}(\Omega)}^{2\gamma-1} \parallel \rho \parallel_{L^{2\beta\gamma+1}(\Omega)}^{\frac{1}{2 \gamma}} \parallel \rho u \parallel_{L^{\frac{2\beta\gamma+1}{\beta}}(\Omega)} \parallel \g u\parallel_{L^2(\Omega)},
\end{eqnarray*}
where we used Proposition~\ref{ME} and H\"{o}lder's inequality.\\
Since
\[\ba
\|\rho u \parallel_{L^{\frac{2\beta\gamma+1}{\beta}}(\Omega)} &\leq C \parallel \rho \parallel_{L^{2\beta\gamma+1}(\Omega)} \parallel u \parallel_{L^{\frac{2\beta\gamma+1}{\beta-1}}(\Omega)}\\
&\leq C\parallel \rho \parallel_{L^{2\beta\gamma+1}(\Omega)}(1+ \parallel \nabla u \parallel_{L^2(\Omega)}),
\ea\]
we obtain that
~\be\label{fgam}
\frac{d}{dt}\int _{\Omega}\rho f^{2\gamma} dx  \leq C\Big(1+ \int_{\Omega}\rho f^{2\gamma } +\int _{\Omega}\rho^{2\beta\gamma+1 }dx\Big) (1+\parallel \g u \parallel_{L^{2}(\Omega)}^2).
\ee

Next, we will estimate the term~$\int_{\Omega}\rho^{2\beta\gamma+1} dx$.\\
Notice that $f=(\theta(\rho) + \xi +G(t))_{+} \geq (\theta(\rho) + \xi +G(t))$ and
$$\rho^{\beta} \leq C\theta(\rho),\text{ when }\rho\geq 2.$$
For any~$\rho\geq 2 $, we have
$$\rho^{\beta}\leq C\theta(\rho) \leq C(f-\xi -G(t)) \leq C(f+|\xi|).$$
Hence, we can devide the term~$\int_{\Omega} \rho^{2\beta\gamma+1} dx$ into two parts:
\begin{eqnarray*}
	\int_{\Omega} \rho^{2\beta\gamma+1} dx &=&\int_{\Omega \cap {\rho (\leq 2)} } \rho^{2\beta\gamma+1} dx +\int_{\Omega \cap (\rho > 2) } \rho^{2\beta\gamma+1} dx\\
	& \leq & C+\int_{\Omega \cap (\rho >2)} \rho f^{2\gamma} dx +C\int_{\Omega}\rho |\xi|^{2\gamma} dx\\
	& \leq & C+\int_{\Omega \cap (\rho >2)} \rho f^{2\gamma} dx  + C\parallel \rho \parallel_{L^{\frac{2\beta\gamma+1}{2\beta\gamma+1-\gamma}}(\Omega)}  \parallel \xi \parallel_{L^{2(2\beta\gamma+1)}(\Omega)}^{2\gamma}\\
	& \leq & C+\int_{\Omega \cap (\rho >2)} \rho f^{2\gamma} dx  + C\parallel \rho \parallel_{L^{\frac{2\beta\gamma+1}{2\beta\gamma+1-\gamma}}(\Omega)}  \parallel \nabla \xi \parallel_{L^{\frac{2\beta\gamma+1}{\beta\gamma+1}}(\Omega)}^{2\gamma}\\
	& \leq & C+\int_{\Omega \cap (\rho >2)} \rho f^{2\gamma} dx  + C\parallel \rho \parallel_{L^{\frac{2\beta\gamma+1}{2\beta\gamma+1-\gamma}}(\Omega)} \parallel\rho^{\frac{1}{2}}\parallel_{L^{2(2\beta\gamma+1)}(\Omega)}^{2\gamma} \parallel\rho ^{\frac{1}{2}} u\parallel_{L^2(\Omega)}^{2\gamma},
\end{eqnarray*}
Since $$\frac{2\beta\gamma+1}{2\beta\gamma+1-\gamma} \leq 2\beta\gamma+1,$$
we have
\begin{eqnarray*}
	\parallel \rho \parallel_{L^{\frac{2\beta\gamma+1}{2\beta\gamma+1-\gamma}}(\Omega)} \parallel\rho^{\frac{1}{2}}\parallel_{L^{2(2\beta\gamma+1)}(\Omega)}^{2\gamma}
	&\leq&\parallel \rho \parallel_{2\beta\gamma+1} \parallel \rho \parallel_{2\beta\gamma+1}^{\gamma}=\parallel \rho \parallel_{2\beta\gamma+1}^{\gamma+1}.
\end{eqnarray*}
For~$\parallel \rho \parallel_{2\beta\gamma+1}^{\gamma+1}$,
taking~$p=\frac{2\beta\gamma+1}{\gamma+1}, \quad q=\frac{2\beta\gamma+1}{2\beta\gamma-\gamma}$ satisfying $\frac{1}{p} + \frac{1}{q}=1$, and using Young's inequality yield
\[\parallel \rho \parallel_{2\beta\gamma+1}^{\gamma+1} \leq \eps (\parallel \rho \parallel_{2\beta\gamma+1}^{\gamma+1})^{\frac{2\beta\gamma+1}{\gamma+1}}+C_\eps.
\]
Take~$\eps$ small enough such that $C\eps\leq \f 12$, we obtain that
\begin{eqnarray*}
	\int_{\Omega} \rho^{2\beta\gamma+1} dx &\leq &C+ \int_{\Omega}\rho f^{2\gamma} dx +C\parallel \rho \parallel_{2\beta\gamma+1}^{\gamma+1}\\
	&\leq & C+ \int_{\Omega}\rho f^{2\gamma} dx +\frac{1}{2}\parallel \rho \parallel_{2\beta\gamma+1}^{2\beta\gamma+1}+C.
\end{eqnarray*}
Therefore,
$$\int_{\Omega} \rho^{2\beta\gamma+1} dx \leq C+C \int_{\Omega} \rho f^{2\gamma} dx.$$
Submitting the above into~\eqref{fgam} gives
~\be
\frac{d}{dt}\int _{\Omega}\rho f^{2\gamma} dx  \leq C\Big(1+ \int_{\Omega}\rho f^{2\gamma }dx \Big) (1+\parallel \div u \parallel_{L^{2}(\Omega)}^2).
\ee
Finally, ~Gronwall's inequality implies that
$$\sup_{t\in[0,\ T]}\int_{\Omega} \rho^{2\beta\gamma+1 }dx \leq \sup_{t\in[0,\ T]}\int _{\Omega}\rho f^{2\gamma} dx +C \leq C.$$
This completes the proof.
\endproof

\begin{Proposition}\label{QS} If
	~\be
	0<\v <\min\Big\{2,\ \f{4\gamma(\beta-1)+2}{2\gamma-1}\Big\},
	\ee
	then there exists a constant $C_\v>0$ depending on the time T such that
	~\be\label{eudelta}
	\sup_{t\in [0,\ T]}\int_{\Omega} \rho |u|^{\v+2} dx\leq C_\v(T).
	\ee
\end{Proposition}

\proof
Multiplying~$\eqref{eq1}_{2}$ by~$(2+\v)u\cdot |u|^{\v}$, and integrating (by parts) over $\Omega$ gives  
~\be\label{Udelta}
\ba
\frac{d}{dt} \int_{\Omega}\rho|u|^{\v+2} dx &+ \mu(\v+2) \int_{\Omega}[\v|u|^{\v}\cdot |\nabla |u||^2 + |u|^{\v} |\nabla u|^{2}]dx\\
&+ (2+\v) \int_{ \Omega }(\mu+ \lambda(\rho)) |u|^{\v} |\div u|^{2}dx\\
&+ (2+\v) \int_{ \Omega }\v|u|^{\v-1}  (u \cdot \nabla) |u| \cdot (\mu+ \lambda(\rho))\div u dx\\
= &\ (2+\v)\int_{\Omega} P \cdot \div (|u|^{\v}\cdot u) dx\\
\leq &\ C\int_{\Omega} \rho^{\gamma} |u|^{\v} |\nabla u| dx.
\ea\ee
Here we used the fact that
\begin{eqnarray*}
	\div (|u|^{\v} \cdot u) &=& u \cdot \nabla |u|^{\v} + |u|^{\v} \cdot \div u\\
	& \leq & C|u|^{\v}  |\nabla u|.
\end{eqnarray*}

Notice that
\[\ba
(2+\v)\v \int_{ \Omega } &(\mu+ \lambda(\rho))\div u |u|^{\v-1}  (u \cdot \nabla) |u|  dx\\
&=(2+\v)\v \int_{ \Omega } (\mu+ \lambda(\rho))|u|^{\v}\big(|\p_r u|^2 + \f u r\p_r u\big) dx\\
&\geq (2+\v)\v \int_{ \Omega } (\mu+ \lambda(\rho))|u|^{\v}\big(|\p_r u|^2 -\f 12|\div u|^2\big) dx,
\ea\]
which combined with~\eqref{Udelta} gives
\[\ba
\frac{d}{dt} \int_{\Omega}\rho|u|^{\v+2} dx &+ \mu(\v+2) \int_{\Omega}[\v|u|^{\v}\cdot |\nabla |u||^2 + |u|^{\v} |\nabla u|^{2}]dx\\
&+ (2+\v) \int_{ \Omega }(\mu+ \lambda(\rho)) |u|^{\v} |\div u|^{2}dx\\
&+ (2+\v)\v \int_{ \Omega } (\mu+ \lambda(\rho))|u|^{\v}\big(|\p_r u|^2 -\f 12|\div u|^2\big) dx\\
\leq & C\int_{\Omega} \rho^{\gamma} |u|^{\v} |\nabla u| dx.
\ea\]
Since~$1-\frac{\v}{2} > 0 $, after applying ~H\"{o}lder's inequality and ~Young's inequality, we obtain that
\[\ba
\frac{d}{dt}\int_{\Omega}\rho |u|^{2+\v} dx &+ \mu(\v+2)\int_{\Omega}|u|^{\v} |\nabla u|^{2} dx
\leq C\int_{\Omega}\rho^{\gamma}|u|^{\v} |\nabla u| dx\\
& \leq \frac{\mu}{2}(\v+2)\int_{\Omega} |u|^{\v}|\nabla u|^{2} dx + C\int_{\Omega} \rho |u|^{2+\v} dx +C\int_{\Omega} \rho^{\gamma(2+\v)-\f \v 2} dx,
\ea\]
Hence, \eqref{eudelta} follows from Proposition~\ref{BJ}, the fact  $\gamma(2+\v)-\f \v 2 \leq 2\beta\gamma+1$, and~Gronwall's inequality immediately.
\endproof

\section{A priori estimates (I): upper bound of the density}

Next, we will give the upper bound of the density.
\begin{Proposition}\label{SZ} Under the conditions of Theorem \ref{them1}, there is a constant $C>0$ depending on the time T, such that
	~\be
	\parallel \rho \parallel_{L^{\infty} L^{\infty}}  \leq C.
	\ee
\end{Proposition}	 

\proof Let ~$\t G(t)\triangleq\int_{0}^{t} F(R, \ts)d\ts$. Recalling the equation of $\theta$,
$$(\theta + \xi + \t G(t))_{t} + u\cdot \nabla(\theta+\xi +\t G(t)) +\int_{R}^{r} \frac{\rho u^2}{s} ds +P(r, t) =0.$$
we have
$$(\theta + \xi + \t G(t))_{t} + u\cdot \nabla(\theta+\xi +\t G(t))  +P(r, t) = -\int_{R}^{r} \frac{\rho u^2}{s} ds
\leq \|\rho\|_{L^\infty L^\infty}\|\g u\|_2^2.
$$
Integrating above inequality on ~$(0,\ t)$, we get
\begin{eqnarray*}
	(\theta + \xi + \t G(t))(r, t) & \leq &\|(\theta + \xi )(r, 0)\|_{\infty}+\|\rho\|_{L^\infty L^\infty}\int_{0}^{t} \|\g u\|_2^2 d\ts\\
	& \leq & C\Big(1 + \parallel \rho \parallel_{L^{\infty} L^{\infty}}\Big).
\end{eqnarray*}
By proposition ~\ref{QS} and ~H\"{o}lder's inequality,
\begin{eqnarray*}
	|\xi| \leq \int_{0}^{R} |\rho u| dr & = & \Big |\int_{0}^{R} \rho^{\frac{1}{2+\v}} \cdot \rho^{\frac{1+\v}{2+\v }} \cdot |u| \cdot r^{\frac{1}{2+\v}} \cdot r^{-\frac{1}{2+\v}} dr \Big|\\
	& \leq &\Big [\int_{0}^{R}(\rho^{\frac{1}{2+\v}} \cdot |u| \cdot r^{\frac{1}{2+\v}} )^{2+\v}dr  \Big]^{\frac{1}{2+\v}} \cdot
	\Big[ \int_{0}^{R} (\rho^{\frac{1+\v}{2+\v}}r^{-\frac{1}{2+\v}} )^{\frac{\v+2}{\v+1}} dr \Big]^{\frac{\v+1}{\v+2}}\\
	&\leq & C(T)\parallel \rho \parallel_{L^{\infty}}^{\frac{1+\v}{2+\v}} (\int_{0}^{R} r^{-\frac{1}{1+\v}} dr)^{\frac{\v+1}{\v+2}}\\
	& \leq &C(T) \parallel \rho \parallel_{L^{\infty}}^{\frac{1+\v}{2+\v}}.
\end{eqnarray*}
Applying lemma ~\ref{FR}, it holds that
\[
|\t G(t)|\leq C(1+ \int_0^t\|\rho\|_\beta^{\beta} d\ts)\leq C(1+ \int_0^t\|\rho\|_\infty^{\beta} d\ts)
\]
Therefore, we get
\begin{eqnarray*}
	\theta &\leq& (\theta+\xi+\t G) +|\xi| + |\t G|\\
	& \leq & C\Big( 1+ \parallel \rho \parallel_{L^{\infty} L^{\infty}} +\parallel \rho \parallel_{L^{\infty}L^{\infty}}^{\frac{1+\v}{2+\v}}+ \int_0^t\|\rho\|_\infty^{\beta} d\ts\Big)\\
	& \leq & C\Big( 1+ \parallel \rho \parallel_{L^{\infty} L^{\infty}}+ \int_0^t\|\rho\|_\infty^{\beta} d\ts\Big).
\end{eqnarray*}

Notice that ~$\theta(\rho)= 2\mu \ln \rho + \frac{1}{\beta}(\rho^{\beta}-1)$, we have
$$\ba
\rho^{\beta} &\leq \max (2^{\beta}, \theta(\rho))\\
&\leq C\Big( 1+ \parallel \rho \parallel_{L^{\infty} L^{\infty}}+ \int_0^t\|\rho\|_\infty^{\beta} d\ts\Big).\\
\ea$$
Therefore, 
~\be\label{SQJ}
\|\rho\|_\infty^{\beta}\leq C\Big( 1+ \parallel \rho \parallel_{L^{\infty} L^{\infty}}+ \int_0^t\|\rho\|_\infty^{\beta} d\ts\Big).
\ee
Applying ~Gronwall's inequality, it follows that
\[
\int_0^t\|\rho\|_\infty^{\beta} d\ts\leq C\Big( 1+ \parallel \rho \parallel_{L^{\infty} L^{\infty}}\Big).
\]
This, together with~\eqref{SQJ}, gives
\[
\|\rho\|_\infty^{\beta}\leq C\Big( 1+ \parallel \rho \parallel_{L^{\infty} L^{\infty}}\Big).
\]
Taking superium according to ~$t\in [0,T]$, one gets
\[
\|\rho\|_{L^{\infty}L^{\infty}}^{\beta}\leq C\Big( 1+ \parallel \rho \parallel_{L^{\infty} L^{\infty}}\Big).
\]
Applying ~Young inequality, we have
\[
\|\rho\|_{L^{\infty}L^{\infty}}^{\beta}\leq \f 12 \|\rho\|_{L^\infty L^\infty}^{\beta} + C.
\]
This gives
$$\parallel \rho \parallel_{L^{\infty}L^{\infty}}^{\beta} \leq C(T),$$
and finishes the proof of Proposition \ref{SZ}.
		  
\endproof

\section{A priori estimates (II): higher order estimates}

Indeed, the bound of the density is enough to extend the local strong solution to a global one. The proof is analogous to the blowup criterion established in \cite{2016Huang}.

To make the content self-contained, we give a detailed proof as follows.
\begin{lem} If~$\|\r\|_\infty \leq C$, then for all~$1\leq q \leq \infty$, it holds that
	~\be\label{ZJ}
	\parallel \nabla F \parallel_{q} \leq C \parallel \nabla \dot{\u} \parallel_{2}.
	\ee
\end{lem}
\proof
Applying Lemma~\ref{ME} to~$$\rho \dot {u}= \nabla F,$$ one has
\[
\parallel \nabla F \parallel_{q} = \parallel \rho \dot{u} \parallel_{q} \leq C\parallel \dot{u} \parallel_{q} \leq C\parallel \nabla \dot{u} \parallel_{2}.
\]
\endproof
Next, we give the estimate of $\parallel \nabla \u \parallel_{4}^{4}$.
\begin{lem}\label{FOUR}
	If~$\|\r\|_\infty \leq C$, then
	~\be
	\| \g \u \| _{4}^{4}\leq C+ C(\|\g \u\|^2 +1)\Big(\|\sqrt{\rho}\dot \u\|_{2}^2 + \int_\O \f{F^2}{2\mu+\ld} dx\Big).
	\ee
	
\end{lem}
\proof 
For $\u(x,t)=u(r,t)\dfrac{x}{r}$, $\g\times \u = 0$, so
\be
\left\{
\ba
&\Delta \u = \nabla \div \u ,\\
& \u |_{\partial \Omega} =0.\\
\ea
\right.
\ee
The standard $L^p$-estimate for the elliptic system yields
\[\ba
\parallel \nabla \u \parallel _{4}^{4} &\leq C \parallel \div \u \parallel_{4}^{4} = C \parallel \frac{F+P}{2 \mu + \lambda} \parallel_{4}^{4}
\leq C\|F\|_4^4\\
& \leq C(1+ \parallel F -\bar F \parallel_{4}^4 + |\bar F|^4)\\
& \leq C(1+ \parallel F -\bar F \parallel_{2}^2 \cdot \parallel \nabla F \parallel_{2}^2 + \|F\|_2^4)\\
& \leq C(1+ \parallel F \parallel_{2}^2 \cdot \parallel\sqrt{\rho}\dot \u\parallel_{2}^2 + \|F\|_2^4)\\
&= C\Big(1+ \| F \|_{2}^2 (\|\sqrt{\rho}\dot \u\|_{2}^2 + \|F\|_2^2)\Big).
\ea\]
Notice that ~$\|F\|_2\leq C(\|\g \u\|_2+1)$ and
$$
\int_\O F^2dx \leq C\int_\O \f{F^2}{2\mu+\ld} dx,
$$
one has
$$\ba
\parallel \nabla \u \parallel _{4}^{4}\leq C\Big(1+ (\|\g \u\|_2^2+1)( \|\sqrt{\rho}\dot \u\|_{2}^2 + \int_\O \f{F^2}{2\mu+\ld} dx)\Big).\\
\ea,$$
which completes the proof of the Lemma.
\endproof

\begin{lem}
	If ~$\|\r\|_\infty \leq C$, then we have
	~\be\label{LL2}
	\ba
	\f 12\f{d}{dt}\int_{\Omega} \f{F^2}{2\mu+\ld}dx &+ \int_{\Omega} \r |\dot{\u}|^2 dx\leq C + C\|\g \u\|_4^4.
	\ea\ee
\end{lem}
\proof		   

Multiplying $\eqref{eq1}_2$ by $\d u$, then integrating the resulting identity over space, one gets after integration by parts and using transport formula
\[\ba
\int_{\Omega} \r |\d u|^2 dx&= \int_{\Omega} \g F\cdot \d u dx = -\int_{\Omega} F\div\d u dx\\
&=-\int_{\Omega} \Big[F\da\div u + F\g u:\g u^t \Big]dx\\
&=-\int_{\Omega} \Big[(2\mu+\ld)\div u \da\div u -P\da\div u + F\g u:\g u^t \Big]dx\\
&=-\f 12\f{d}{dt}\int_{\Omega} \Big[(2\mu+\ld)\big|\div u -\f{P}{2\mu+\ld}\big|^2 - \f{P^2}{2\mu+\ld} \Big]dx \\
&\quad-\int_{\Omega} \Big[\r|\div u|^2\big(\f{2\mu+\ld}{2\r}\big)_\r\r\div u  - \r\div u \big(\f{P}{\r}\big)_\r\r\div u + F\g u:\g u^t \Big]dx.
\ea\]
Since
\[\ba
\frac{1}{2}\frac{d}{dt} \int_{\Omega} \frac{P^2}{ 2 \mu+ \lambda} dx&= \frac{1}{2} \int_{\Omega} \frac{\partial}{\partial t} (\frac{P^2}{ 2 \mu+ \lambda}) dx
= \frac{1}{2}\int_{\Omega}(\frac{P^2}{ 2 \mu+ \lambda})_{\rho} \cdot \frac{D}{Dt} \rho dx\\
&=-\frac{1}{2} \int_{\Omega}(\frac{P^2}{ 2 \mu+ \lambda})_{\rho} \cdot \rho \div u dx,
\ea\]
we have
\[\ba
\int_{\Omega} \rho |\dot{u}|^2 dx &=-\f 12\f{d}{dt}\int_{\Omega} (2\mu+\ld)\big|\div u  -\f{P}{2\mu+\ld}\big|^2  dx-\f 12\int_{\Omega} \big(\f{P^2}{2\mu+\ld}\big)_\r\r\div udx \\
&\quad+\int_{\Omega} \Big[-\r|\div u|^2\big(\f{2\mu+\ld}{2\r}\big)_\r\r\div u  - \r\div u \big(\f{P}{\r}\big)_\r\r\div u + F\g u:\g u^t \Big]dx\\
&=-\f 12\f{d}{dt}\int_{\Omega} \f{F^2}{2\mu+\ld}dx -\f 12\int_{\Omega} \big(\f{P^2}{2\mu+\ld}\big)_\r\r\div udx \\
&\quad+\int_{\Omega}-\r|\div u|^2\big(\f{2\mu+\ld}{2\r}\big)_\r\r\div u  - \r\div u \big(\f{P}{\r}\big)_\r\r\div u + F\g u:\g u^t dx\\
&\leq -\f 12\f{d}{dt}\int_{\Omega} \f{F^2}{2\mu+\ld}dx + C\|\div u\|_2^2 + C\int_{\Omega}|\div u|^3 dx + \int_{\Omega} |F||\g u|^2 dx
\ea\]

Hence, H$\ddot{o}$lder's inequality leads to
\[\ba
\f 12\f{d}{dt}\int_{\Omega} \f{F^2}{2\mu+\ld}dx &+ \int_{\Omega} \r |\d u|^2\leq C\|\div u\|_2^2 + C\int_{\Omega}|\div u|^3 dx + \int_{\Omega} |F||\g u|^2 dx\\
&\leq C(1 + \|\g u\|_4^4).
\ea\]
Hence, we finishes the proof.
\endproof 

\begin{lem} \label{JS}
	If $\|\r\|_\infty \leq C$, it holds that
	~\be\label{LL1}
	\ba
	\f 12\f{d}{dt}\|\sqrt{\r}\d u\|_2^2 &+ \mu\|\div\d u\|_2^2 + \f 12\int_{\Omega}\ld\Big|\da{\div u}\big|^2dx \leq  C\|\g u\|_4^4 + C.
	\ea\ee
\end{lem}	
\proof 
Rewrite the momentum equation $\eqref{eq1}_2$ as
~\be\label{MM}
\rho \dot{u} + \nabla P = L u,
\ee
which yields
\[
\dot{ \rho} \dot{u} + \rho \ddot{u} + \dot{\nabla P} = \dot{Lu}.
\]
So we have
\[\ba
\rho \ddot{u} &= \dot{Lu} - \dot{ \rho} \dot{u}- \dot{\nabla P} \\
&=\dot{Lu} + \rho  \div u \cdot\dot{u}- \dot{\nabla P} \\
&=(Lu)_{t} + \div (Lu \otimes u) - \nabla P_{t} - \div (\nabla P \otimes u).\label{QY}
\ea\]
Multiplying the above equation by $\dot{u}$ and integrating the resulting equality over $\O,$ one gets
~\be\label{WP}
\int_{\Omega} \rho \ddot{u} \dot{u} dx = \int_{\Omega} [(Lu)_{t}+\div (Lu \otimes u)]\cdot \dot u dx - \int_{\Omega} [\nabla P_{t} + \div (\nabla P \otimes u)]\cdot \dot u  dx.
\ee
For the second term of the right hand side of \eqref{WP}, after using $\eqref{eq1}_1$ and integration by parts, we have
	
\[\ba
-\int_{\Omega}(\nabla P_{t} &+ \div (\nabla P \otimes u )) \cdot \dot{u} dx = \int_{\Omega} P_{t} \div \dot{u} + \nabla P (u \cdot \nabla) \dot{u} dx\\
&= \int_{\Omega} P^{'} \rho_{t} \div \dot{u} + \nabla P (u \cdot \nabla) \dot{u} dx\\
&=-\int_{\Omega} \rho P^{'} \div u \div \dot {u} dx- \int_{\Omega} P[\partial_{k}u_{j} \partial_{j} \dot{u_{k}} - \partial_{j} u_{j} \partial_{k} \dot{u_{k}} ]dx,
\ea\]
which leads to
~\be\label{PT2}
\ba
\Big| \int_{\Omega} (\nabla P_{t} + \div (\nabla P \otimes u)) \cdot \dot{u} dx\Big|
& \leq C \parallel \nabla u \parallel _{2} \parallel \nabla \dot{u} \parallel_{2} = C \parallel \nabla u \parallel _{2} \parallel \div  \dot{u} \parallel_{2}.\\
&\leq \frac{\mu}{8} \parallel \div \dot{u} \parallel _{2}^2 + C\parallel \nabla u \parallel_{2}^{2}.
\ea\ee
For the first term of the right hand side of \eqref{WP},
~\be\label{WPr12}
\ba
\int_{\Omega} [(Lu)_{t} &+\div (Lu \otimes u)]\cdot \dot u dx
=\int_{\Omega}(Lu)_{t}\cdot{\dot u} dx + \int_{\Omega} \div (Lu \otimes u) \cdot{\dot u}dx\\
&= 2 \mu \int_{\Omega} \Big[(\nabla \div u)_{t} + \div (\nabla \div u \otimes u)\Big] \cdot \dot{u} dx \\
&\qquad+ \int_{\Omega} \Big[(\nabla( \lambda \div u))_{t} + \div (\nabla(\lambda \div u) \otimes u)\Big] \cdot \dot{u} dx.
\ea\ee
For the first term of the right hand side of \eqref{WPr12}, notice that 
\[\ba
\nabla(u\div u ) 
&= \nabla \div u \otimes u + \div u \nabla u,
\ea\]
\[\ba
\int_{\Omega} \div (\nabla(u \div u)) \cdot \dot{u} dx &= \int_{\Omega} \nabla (\div(u \div u)) \cdot \dot{u} dx\\
&= - \int_{\Omega} \div(u \div u)\cdot \div \dot{u} dx,
\ea\]
and 
\[\ba
-\int_{\Omega} \div(\div u \nabla u)\cdot \dot{u} dx 
&=\int_{\Omega} \div u \nabla u:\nabla u^{t}dx.
\ea\]
Using H$\ddot{o}$lder's inequality, we obtain that 
\be
\ba\label{CD}
2\mu\int_{\Omega}\Big[(\g\div u)_t &+ \div(\g(\div u)\otimes u)\Big]\cdot\d u dx\\
&= 2\mu\int_{\Omega} -\div u_t\div\d{u} - \div(u\div u)\div\d u + \div u \g u:\g\d u^tdx\\
&\leq 2\mu\int_{\Omega}\Big[ -\div u_t\div\d{u} - \div(u\div u)\div\d u \Big]dx + C\|\g u\|_4^2\|\g\d u\|_2\\
&\leq 2\mu\int_{\Omega} \Big[-\div u_t\div\d{u} - (u\cdot\g)\div u\div\d u - |\div u|^2\div\d u \Big]dx+ C\|\g u\|_4^2\|\div \d u\|_2\\
&\leq 2\mu\int_{\Omega} \Big[-\div u_t\div\d{u} - \div (u\cdot\g u)\div\d u + \g u:\g u^t\div \d u \Big] dx + C\|\g u\|_4^2\|\div\d u\|_2\\
&\leq 2\mu\int_{\Omega} -|\div\d u|^2 dx + C\|\g u\|_4^2\|\div\d u\|_2\\
&\leq -\f{3}{2}\mu\|\div\d u\|_2^2 + C\|\g u\|_4^4,
\ea\ee
where we used the fact that
$$
\div (u \cdot \nabla u) = \nabla u : \nabla u^{t} + (u \cdot \nabla) \div u.
$$
	  

Now we come to estimate the second term of the right hand side of \eqref{WPr12}. Since
\[\ba
\int_{\Omega} \nabla (\lambda \div u)_{t} \cdot \dot{u} dx 
&= -\int_{\Omega} (\lambda \div u)_{t} \cdot \div \dot{u} dx,
\ea\]
\[\ba
\nabla (u \lambda \div u) 
&= \nabla (\lambda \div u) \otimes u + \lambda \div u \nabla u,
\ea\]
and 
\[\ba
-\int_{\Omega} \div (\lambda \div u \cdot \nabla u) \cdot \dot{u} dx 
&=\int_{\Omega} \lambda \div u \nabla u : \nabla \dot{u}^t dx,
\ea\]
we get
\[\ba
\int_{\Omega}\Big[&(\g(\ld\div u))_t + \div(\g(\ld\div u)\otimes u)\Big]\cdot\d u dx\\
&= \int_{\Omega} \Big[-(\ld\div u)_t\div\d{u} - \div(u\ld\div u)\div\d u + \ld\div u \g u:\g\d u^t \Big]dx\\
&\leq \int_{\Omega} \Big[-\ld\div u_t\div\d{u} - \div(u\ld\div u)\div\d u -\ld_t\div u\div\d u \Big]dx + C\|\g u\|_4^2\|\g\d u\|_2\\
&= \int_{\Omega} \Big[-\ld\div u_t\div\d{u} - u\cdot\g(\ld\div u) \div\d u - \ld|\div u|^2\div \d u\\
&\qquad -\ld_t\div u\div\d u \Big] dx + C\|\g u\|_4^2\|\div\d u\|_2\\
&\leq \int_{\Omega} \Big[-\ld\div u_t\div\d{u} - \ld u\cdot\g\div u \div\d u - u\cdot\g\ld\div u\div\d u\\
&\qquad -\ld_t\div u\div\d u\Big] dx + C\|\g u\|_4^2\|\div\d u\|_2\\
&\leq \int_{\Omega} -\ld\da{\div u}\div\d{u} -(\ld_t + u\cdot\g\ld)\div u\div\d u  dx + C\|\g u\|_4^2\|\div\d u\|_2.
\ea\]
Observing that
\[\ba
\lambda_{t} + u \cdot \nabla \lambda 
&= -\lambda^{'} \rho \div u,
\ea\]
we have
\[\ba
\int_{\Omega}\Big[(\g(\ld\div u))_t &+ \div(\g(\ld\div u)\otimes u)\Big]\cdot\d u dx\\
&\leq \int_{\Omega} -\ld\da{\div u}\div\d{u} +\r\ld'|\div u|^2\div\d u  dx + C\|\g u\|_4^2\|\div\d u\|_2\\
&\leq \int_{\Omega} -\ld\da{\div u}\div\d{u} dx + C\|\g u\|_4^2\|\div\d u\|_2\\
&\leq \int_{\Omega} -\ld\da{\div u}\div\d{u} dx + \f{\mu}{8}\|\div\d u\|_2^2 +  C\|\g u\|_4^4.
\ea\]

For the first term of the above inequality, 
\[\ba
\int_{\Omega} \ld\da{\div u}\div\d{u} dx &= \ld\da{\div u}\Big[\div u_t + \div(u\cdot\g u) \Big]dx\\
&= \int_{\Omega}\ld\da{\div u}\Big[\div u_t + (u\cdot\g)\div u + \g u:\g u^t\Big]dx\\
&= \int_{\Omega}\ld\Big|\da{\div u}\big|^2 + \ld\da{\div u}\g u:\g u^tdx\\
&\geq \f 12\int\ld\Big|\da{\div u}\big|^2dx - C\|\g u\|_4^4,
\ea\]
where we used the fact that
\[
\lambda a^2 + \lambda a b \geq \frac{\lambda}{2} a^2 - C b^2.
\]
Hence we have
~\be\ba\label{HYL}
\int_{\Omega}\Big[(\g(\ld\div u))_t &+ \div(\g(\ld\div u)\otimes u)\Big]\cdot\d u dx\\
&\leq -\f 12\int_{\Omega}\ld\Big|\da{\div u}\big|^2dx + C\|\g u\|_4^2+\frac{\mu}{8}\|\div\d u\|_2 + C\|\g u\|_4^4\\
&\leq -\f 12\int_{\Omega}\ld\Big|\da{\div u}\big|^2dx + \f{\mu}{8}\|\div\d u\|_2^2 + C\|\g u\|_4^4.
\ea\ee

Combining~\eqref{CD} and \eqref{HYL} yields
\be\label{TKX}
\ba
\int_{\Omega} \Big[(Lu)_t &+ \div(Lu\otimes u)\big]\cdot\d u dx
\leq -\frac{3}{2} \mu \parallel \div \dot{u} \parallel_{2}^{2} + C\parallel \nabla u \parallel_{4}^{4} \\
& \qquad- \frac{1}{2} \int_{\Omega} \lambda \Big| \frac{D}{Dt} \div u \Big|^2dx
 + \frac{\mu}{8} \parallel \div \dot{u}\parallel_{2}^2 + C\parallel \nabla u \parallel_{4}^{4}\\
&\leq -\f{5}{4}\mu\|\div\d u\|_2^2-\f 12\int_{\Omega}\ld\Big|\da{\div u}\big|^2dx + C\|\g u\|_4^4.
\ea\ee
It follows from transport formula \eqref{CJ} that
~\be\label{WZDS}
\int_{\Omega} \rho \ddot{u} \dot{u} dx = \frac{1}{2} \int_{\Omega}\rho \frac{D}{Dt}|\dot{u}|^2 dx = \frac{1}{2} \frac{d}{dt} \int_{\Omega} \rho |\dot{u}|^2dx.
\ee
Inserting~\eqref{PT2}, \eqref{WZDS} and \eqref{TKX} into \eqref{WP} leads to
~\be\ba
\frac{1}{2} \frac{d}{dt} \int_{\Omega} \rho |\dot{u}|^2 dx &\leq -\frac{5}{4} \mu \parallel \div \dot{u}\parallel_{2}^{2} - \frac{1}{2} \int_{\Omega} \lambda \Big|\frac{D}{Dt} \div u \Big|^2 dx + C \parallel \nabla u \parallel_{4}^{4}\\
&+ \frac{\mu}{8} \parallel \div \d u \parallel_{2}^{2} + C\parallel \nabla u \parallel_{2}^{2},
\ea\ee
which implies \eqref{LL1}.
\endproof
\begin{Proposition} If ~$\|\r\|_\infty \leq C$, then we have
	~\be\label{LL3}
	\|\sqrt{\r}\d \u\|_{L^\infty L^2}\leq C,\quad \|F\|_{L^\infty L^2}\leq C
	\ee
	and
	~\be\label{LL4}
	\|\g\dot \u\|_{L^2L^2}\leq C.
	\ee
\end{Proposition}
\proof 
Combining \eqref{LL1}, \eqref{LL2} and Lemma \ref{FOUR} yields
\[\ba
\f 12\f{d}{dt}\Big(\|\sqrt{\r}\d u\|_2^2 &+ \int_{\Omega} \f{F^2}{2\mu+\ld}dx \Big) + \mu\|\div\d u\|_2^2 + \f 12\int_{\Omega}\ld\Big|\da{\div u}\big|^2dx
+\int_{\Omega} \r |\d u|^2dx\\
&\leq C + C\|\g u\|_4^4\\
&\leq  C + C(\|\g u\|^2+1) \Big(\|\sqrt{\rho}\dot u\|_{2}^2 + \int_\O \f{F^2}{2\mu+\ld} dx\Big).
\ea\]
Applying Gronwall's inequality  and \eqref{HD} leads to \eqref{LL3} ~\eqref{LL4}.
\endproof

We still need the following lemma concerning the higher order estimates on $(\rho, \u).$

\begin{lem}\label{lem5} If~$\|\r\|_\infty \leq C$, it holds that
	~\be\label{ru1}
	\|\g\r\|_q \leq C,\quad\text{ and }\quad \|\g^2 u\|_2 \leq C.
	\ee
\end{lem}
\proof	   
It follows from $\eqref{eq1}_{1}$ that
~\be
(\g\r)_t + (u\cdot\g)\g\r + \g u\cdot\g \r + \g\r\div u + \rho\g\div u= 0.
\ee
Multiplying the above equation by $q|\g\r|^{q-2}\g\r$, and integrating the resulting equation over $\O$, we have
\[\ba
\frac{d}{dt} \int_{\Omega} |\nabla \rho|^{q} dx &\leq (q-1)\int_{\Omega} |\nabla \rho|^{q} |\div u|dx + \int_{\Omega} |\nabla \rho|^{q} |\nabla u| dx + \int_{\Omega} |\nabla \rho|^{q-1} \cdot |\nabla \div u| dx\\
&\leq C\int_{\Omega}|\nabla \rho|^{q} |\nabla u| dx + C\int_{\Omega}|\nabla \rho|^{q-1} |\nabla \div u| dx\\
&\leq C \|\g u\|_\infty \|\g\r\|_q^q + C \|\g\div u\|_q\|\g\r\|_q^{q-1},
\ea\]
which implies that
~\be\label{f54}\ba
\f{d}{dt}\|\g\r\|_q & \leq C \|\g u\|_\infty \|\g\r\|_q + C \|\g\div u\|_q\\
&\leq  C \big[\|\div u\|_\infty \log(e + \|\g^2 u\|_q) + \|\g u\|_2 + 1\big] \|\g\r\|_q + C \|\g\div u\|_q.
\ea\ee
From the definition of the effective viscous flux $F= (2 \mu+ \lambda)\div u -P$, we have
~\be\label{divi}\ba
\|\div u\|_\infty & \leq \|\f{F+p}{2\mu + \ld}\|_\infty
\leq C(\|F\|_\infty + 1)\\
&\leq C( \|F\|_2^{\f{q-2}{2q-2}}\|\g F\|_q^{\f{q}{2q-2}} +\|F\|_2+1)\\
&\leq C( \|\g F\|_q^{\f{q}{2q-2}} +1)\\
&\leq C( \|\r\d u\|_q^{\f{q}{2q-2}} +1)\\
&\leq C( \|\g\d u\|_2^{\f{q}{2q-2}} +1),
\ea\ee
where we used the fact that
\be
\parallel \nabla F \parallel_{q} \leq C\parallel \rho \dot{u} \parallel_{q} \leq C\parallel \dot{u} \parallel_{q} \leq C\parallel \nabla \dot{u} \parallel_{2}.
\ee
Next, since 
~\be
\left\{
\ba
\Delta u = \nabla \div u ,\\
u |_{\partial \Omega} =0,\\
\ea
\right.
\ee
by the~$L^p$ theory for the elliptic system (see~\cite{GT}), we obtain that
\[\ba
\|\g^2 u\|_q &\leq C\|\g\div u\|_q
\leq C\|\g\big(\f{F+p}{2\mu + \ld}\big)\|_q\\
&\leq C(\|\g F\|_q + \|\g\r\|_q + \|F\|_\infty\|\g\r\|_q)\\
&\leq C(\|\g\d u\|_2 + \|\g\r\|_q + \|F\|_2^{\f{q-2}{2q-2}}\|\g F\|_q^{\f{q}{2q-2}}\|\g\r\|_q)\\
&\leq C(\|\g\d u\|_2 + \|\g\r\|_q + \|\g\d u\|_2^{\f{q}{2q-2}}\|\g\r\|_q)\\
&\leq C(\|\g\d u\|_2 + \|\g\r\|_q + \|\g\r\|_{q}^{\f{2q-2}{q-2}})\\
&\leq C(e + \|\g\d u\|_2 + \|\g\r\|_q^{\f{2q-2}{q-2}}).
\ea\]
Here we used the facts that
\[\ba
\parallel  \nabla \dot{u} \parallel_{2}^{\frac{q}{2q-2}} \cdot \parallel \nabla \rho \parallel_{q}
&\leq \parallel  \nabla \dot{u} \parallel_{2}^{\frac{q}{2q-2} \cdot \frac{2q-2}{q}} + \parallel \nabla \rho \parallel_{q}^{\frac{2q-2}{q-2}}\\
&=\parallel \nabla \dot{u} \parallel_{2} + \parallel \nabla \rho \parallel^{\frac{2q-2}{q-2}}.
\ea\]
and 
\[
\parallel \nabla \rho \parallel_{q} \leq C(e+ \parallel \nabla \rho \parallel_{q}^{\frac{2q-2}{q-2}}), ~\text{for}~\frac{2q-2}{q-2}\geq 1.
\]
Notice that 
\[
\log(e+a+b^\a) \leq C_\a\log(e+a+b), \text{for all} ~a,b>0,\a>1.
\]
we get
\[\ba
\log(e + \|\g^2 u\|_q) &\leq C\log(e + \|\g\d u\|_2 + \|\g\r\|_q^{\f{2q-2}{q-2}})\\
&\leq C\log(e + \|\g\d u\|_2+\|\g\r\|_q).
\ea\]
From~\eqref{divi}, we have~$$\parallel \div u \parallel_{\infty} \leq C(\parallel \rho \dot {u}\parallel_{q}^{\frac{q}{2q-2}}+1),$$
and hence
\[\ba
\|\g\div u\|_q &=\|\g\big(\f{F+p}{2\mu + \ld}\big)\|_q\\
&\leq C(\|\g F\|_q + \|\g\r\|_q + \|F\|_\infty\|\g\r\|_q)\\
&\leq C(\|\g F\|_q + \|\g\r\|_q + \| \div u\|_\infty\|\g\r\|_q)\\
&\leq C(\|\g F\|_q + \|\g\r\|_q + (1+ \parallel \rho \dot{u}\parallel_{q}^{\frac{q}{2q-2}})\|\g\r\|_q)\\
&\leq C(\|\g\d u\|_2 + \|\g\r\|_q + \|\g\d u\|_2\|\g\r\|_q)\\
&\leq C\big[\|\g\d u\|_2 + (1+ \|\g\d u\|_2)\|\g\r\|_q\big].
\ea\]
Inserting the above estimates into~\eqref{f54}, we get
~\be\label{dr1}\ba
&\f{d}{dt}\|\g\r\|_q \leq  C \big[( \|\g\d u\|_2^{\f{q}{2q-2}} +1) \log(e + \|\g\d u\|_2 + \|\g\r\|_q) + \|\g u\|_2 + 1\big] \|\g\r\|_q \\
&\qquad\qquad + C\big[\|\g\d u\|_2 + (1+ \|\g\d u\|_2)\|\g\r\|_q\big]\\
&\leq  C \big[( \|\g\d u\|_2 +1) \log(e + \|\g\r\|_q) + \|\g u\|_2 + \|\g \d u\|_2 + 1\big] \|\g\r\|_q + C\|\g\d u\|_2.
\ea\ee
Set
\[
y(t) = e + \|\g\r\|_q \geq e,
\]
then we have
\[
\log y(t) \geq 1.
\]
It follows from \eqref{dr1} that
\[\ba
\frac{d}{dt}(\parallel \nabla \rho \parallel_{q} &+ e) 
\leq C \big[(\parallel \nabla \dot{u}\parallel_{2} +1) \log(e + \parallel \nabla \rho \parallel_{q}) \\
&+ \parallel \nabla u\parallel_{2} + \parallel \nabla \dot{u} \parallel_{2} +1\big] (\parallel \nabla \rho \parallel_{q} + e)+ C\|\g\d u\|_2,
\ea\]
that is
\[\ba
\f{d y(t)}{dt}& &\leq  C \big[( \|\g\d u\|_2 +1) \log y(t) + \|\g u\|_2 + \|\g \d u\|_2 + 1\big] y(t) + C\|\g\d u\|_2.
\ea
\]
Multiplying the above inequality by $\dfrac{1}{y(t) \log y(t)}$ gives
\[\ba
\frac{d}{dt}(\log \log y(t))&=\frac{1}{y(t) \log y(t)} \frac{d}{dt}  y(t) \\
&\leq C \Big[(\parallel \nabla \dot{u}\parallel_{2} +1)\log y(t) +\parallel \nabla u \parallel_{2} + \parallel \nabla \dot{u} \parallel_{2}+1  \Big] \frac{1}{\log y(t)} \\
&\qquad\qquad + C \parallel \nabla \dot {u} \parallel_{2} \frac{1}{y(t) \log y(t)}\\
&\leq C \Big[(\parallel \nabla \dot{u} \parallel_{2} +1) + \parallel \nabla u \parallel_{2} + \parallel \nabla \dot{u} \parallel_{2} +1 \Big] + C \parallel \nabla \dot{u} \parallel_{2}.
\ea\]
Integrating over~$(0,T)$ yields
\[
\log\log y(T) \leq \log\log y(0) + C,
\]
which implies that
\[
y(T) \leq C.
\]
and hence ~$\parallel \nabla \rho \parallel_{q} \leq C.$

Now we turn to estimate~$\g^2 u$.
By~Sobolev's inequality~\eqref{Sob}, we have
\[\ba
\|\g^2 u\|_2&\leq C\|\g\div u\|_2\leq C\|\g\big(\f{F+p}{2\mu+\ld}\big)\|_2\\
&\leq C(\|\g\r\|_2 + \|\g F\|_2 +  \|F\|_{\f{2q}{q-2}}\|\g \r\|_q)\\
&\leq C\big[\|\g\r\|_2 + \|\g F\|_2 +  (\|F\|_2+ \|\g F\|_2) \|\g \r\|_q\big]\\
&\leq C(\|\g F\|_2+1) \leq C(\|\r \d u\|_2+1) \leq C(\|\sqrt{\r}\d u\|_2+1)\leq C.
\ea\]
\endproof
This is enough to gurantee the strong solution to be a global one.

In the next section, we will show how to derive the uniform boundness of the density independent of the time as long as  
\be
\max\{1,\frac{\gamma+2}{4}\}<\beta\le\gamma.
\ee

\section{Upper bound of the density independent of time}
In this section, we denote $R_T=1+\sup\limits_{0<t<T}\r$. 

Here we will borrow the steps adapted in the proof of \cite{2016Huang} and carefully deal with extra terms caused by the commutator's estimate \eqref{thxi}.

The first Lemma is the corresponding version of Lemma 3.1 in  \cite{2016Huang} which can be much simplied for the radially symmetric case.

\begin{lemma} {\rm (Estimate of $F$)} For any $\alpha\in (0,1)$, there is a constant $C(\alpha)$ depending only on $\alpha, \mu,\beta,\gamma, \|\rho_0\|_{L^\infty}$, and $\|u_0\|_{H^1}$ such that
	\be\label{A2}
	\sup\limits_{0\leq t\leq T}\log\Big(e + \int\f{F^2}{2\mu+\ld} dx + \int(2\mu+\ld)|\div u|^2dx\Big) \leqs C(\a)R_T^{\a\beta+1} + R_T^{\gamma-2\beta+ 1}.
	\ee
\end{lemma}
\proof
The momentum equations can be rewritten as
\[
\r\dot u = \g F.
\]
Consequently, mutliplying the above equations by $\dot{u}$ and integrating the resulting equality over the ball $\Omega$ yields that
\[
\int \r |\dot u|^2 dx = \int \g F\cdot \dot u dx = -\int F\div \dot u dx.
\]
Note that
\[
\div \dot u = \f{D}{Dt} \div u + |\div u|^2 - 2\p_r u\f{u}{r}.
\]
Hence,
\[
\int F\div \dot u dx = \int F \f{D}{Dt} \div udx + \int F|\div u|^2dx - 2\int F\p_r u\f{u}{r}dx
\]
The first term can be estimated as follows
\[\ba
\int &F \f{D}{Dt} \div u dx  =\int F\f{D}{Dt} (\f{F+P}{2\mu+\ld})dx\\
&=\f 12\f{d}{dt}\int\f{F^2}{2\mu+\ld} dx -\f 12\int \f{F^2}{2\mu+\ld}\div u dx\\
&\qquad + \f 12\int F^2 \f{\r\ld'}{(2\mu+\ld)^2}\div udx  - \int F  \r(\f{P}{2\mu+\ld})'\div udx\\
\ea\]
Therefore,
\[\ba
\int \r |\dot u|^2 dx & = -\int F\div \dot u dx\\
&=-\int F \f{D}{Dt} \div udx - \int F|\div u|^2dx + 2\int F\p_r u\f{u}{r}dx\\
&=-\f 12\f{d}{dt}\int\f{F^2}{2\mu+\ld} dx +\f 12\int \f{F^2}{2\mu+\ld}\div u dx- \f 12\int F^2 \f{\r\ld'}{(2\mu+\ld)^2}\div udx \\
&\qquad  + \int F  \r(\f{P}{2\mu+\ld})'\div udx
- \int F|\div u|^2dx + 2\int F\p_r u\f{u}{r}dx\\
\ea\]
Combing all the above equalities together to get
\[\ba
\f 12\f{d}{dt}\int\f{F^2}{2\mu+\ld} dx &+ \int \r |\dot u|^2 dx = \f 12\int \f{F^2}{2\mu+\ld}\div u dx - \f 12\int F^2 \f{\r\ld'}{(2\mu+\ld)^2}\div udx \\
&\qquad + \int F  \r(\f{P}{2\mu+\ld})'\div udx
- \int F|\div u|^2dx + 2\int F\p_r u\f{u}{r}dx\\
&\leqs \int  \f{F^2}{(2\mu+\ld)}|\div u|dx  + \int (\f{|F|P}{2\mu+\ld})|\div u|dx \\
&\qquad + \int F|\div u|^2dx + \|\g F\|_2\|\div u\|_2^2\\
&\leqs \int  \f{F^2}{(2\mu+\ld)}|\div u|dx  + \int (\f{|F|P}{2\mu+\ld})|\div u|dx  + \|\g F\|_2\|\div u\|_2^2\\
&= I_1 + I_2 + I_3.
\ea\]

Then each $I_i$ can estimated as follows:

For any $0<\a<1$,
\[\ba
I_1 &=\int  \f{F^2}{(2\mu+\ld)}|\div u|dx \leqs \|\f{F^2}{2\mu+\ld}\|_2\|\div u\|_2\\
&\leqs \|\f{F}{\sqrt{2\mu+\ld}}\|_2^{1-\a} \|F\|^{1+\a}_{\f{2(1+\a)}{\a}}\|\div u\|_2\\
&\leqs \|\f{F}{\sqrt{2\mu+\ld}}\|_2^{1-\a} \|F\|_2^\a \|\g F\|_2\|\div u\|_2\\
&\leqs R_T^{\f {\a\beta}2}\|\f{F}{\sqrt{2\mu+\ld}}\|_2 \|\g F\|_2 \|\div u\|_2\\
&\leqs R_T^{\f {\a\beta}2}\|\f{F}{\sqrt{2\mu+\ld}}\|_2\|\rho \dot u\|_2\|\div u\|_2\\
&\leq \f 14\|\sqrt{\r} \dot u\|_2^2 +  C(\a)R_T^{\a\beta+1}\|\f{F}{\sqrt{2\mu+\ld}}\|_2^2\|\div u\|_2^2.
\ea\]
Then we will deal with $I_2$.
\[\ba
I_2 &=\int \f{|F|P}{2\mu+\ld}|\div u|dx \leq \|F\|_{2+4\gamma/\beta}\|\f{P}{(2\mu+\ld)^{\f 32}}\|_{2+\beta/\gamma}\|\sqrt{\mu+\ld}\div u\|_2\\
&\leqs R_T^{\gamma/2-\beta}\|\g F\|_2\|\sqrt{\mu+\ld}\div u\|_2\\
&\leqs R_T^{\gamma/2-\beta+ 1/2}\|\sqrt{\r} \dot u\|_2\|\sqrt{\mu+\ld}\div u\|_2\\
&\leq \f 14 \|\sqrt{\r} \dot u\|_2^2 + CR_T^{\gamma-2\beta+ 1}\|\sqrt{\mu+\ld}\div u\|_2^2.\\
\ea\]\\
Then we are about to estimtate $I_3$.
\[\ba
I_3 &= \|\g F\|_2\|\div u\|_2^2 =\|\r \dot u\|_2\|\f{F+P}{2\mu+\ld}\|_2\|\div u\|_2\\
&\leq \f 14 \|\sqrt{\r} \dot u\|_2^2 + CR_T\|\f{F+P}{2\mu+\ld}\|_2^2\|\div u\|_2^2\\
&\leq \f 14 \|\sqrt{\r} \dot u\|_2^2 + CR_T\Big(\|\f{F}{2\mu+\ld}\|_2^2+ \|\f{P}{2\mu+\ld}\|_2^2\Big)\|\div u\|_2^2\\
&\leq \f 14 \|\sqrt{\r} \dot u\|_2^2 + CR_T\|\f{F}{\sqrt{2\mu+\ld}}\|_2^2\|\div u\|_2^2+ CR_T\|\f{P}{2\mu+\ld}\|_2^2\|\div u\|_2^2\\
&\leq \f 14 \|\sqrt{\r} \dot u\|_2^2 + CR_T\|\f{F}{\sqrt{2\mu+\ld}}\|_2^2\|\div u\|_2^2+ CR_T^{\gamma-2\beta+1}\|\div u\|_2^2
\ea\]

Then,
\[\ba
\f{d}{dt}\int\f{F^2}{2\mu+\ld} dx &+ \int \r |\dot u|^2 dx
\leqs   C(\a)R_T^{\a\beta+1}\|\f{F}{\sqrt{2\mu+\ld}}\|_2^2\|\div u\|_2^2\\
&\qquad + CR_T\|\f{F}{\sqrt{2\mu+\ld}}\|_2^2\|\div u\|_2^2 + CR_T^{\gamma-2\beta+1}\|\div u\|_2^2\\
&\qquad + CR_T^{\gamma-2\beta+ 1}\|\sqrt{\mu+\ld}\div u\|_2^2
\ea\]

By Gronwall's inequality, we finally arrive at
\[
\sup\limits_{0\leq t\leq T}\log(e + \int\f{F^2}{2\mu+\ld} dx) \leqs C(\a)R_T^{\a\beta+1} + R_T^{\gamma-2\beta+ 1}.
\]
and moreover,
\[\ba
\int (2\mu+\ld)|\div u|^2 dx &=\int (2\mu+\ld)|\f{F+P}{2\mu+\ld}|^2 dx\\
&\leqs \int \f{F^2}{2\mu+\ld}dx + \int \f{P^2}{2\mu+\ld}dx\\
&\leqs \int \f{F^2}{2\mu+\ld}dx + R_T^{\gamma-\beta}\\
\ea\]
This gives the desired estimates of \eqref{A2} and thus finished the proof Lemma 6.1.
\endproof

Another crucial estimate lies on the $\xi$.
\begin{lemma} {\rm (Estimate of $\xi$)} For any $\alpha\in (0,1)$, there is a constant $C(\alpha)$ depending only on $\alpha, \mu,\beta,\gamma, \|\rho_0\|_{L^\infty}$, and $\|u_0\|_{H^1}$ such that
\[
\|\xi\|_\infty\le C(\a)R_T^{1+\a\beta} + R_T^{1+\gamma/2-\beta}.
\]
\end{lemma}

\proof
Indeed, by Lemma 2.5 , we obtain that for any $\alpha\in (0,1)$
	\[\ba
	\|\xi\|_\infty &\leqs \|\g\xi\|_2\log(e+\|\g\xi\|_3) + \|\xi\|_2 +1\\
	&\leqs \|\r u\|_2\log(e+\|\r u\|_3) + \|\r u\|_{\f{2\gamma}{\gamma+1}} +1\\
	&\leqs R_T^{1/2}\log(e+R_T\|\g u\|_2)  +1\\
	&\leqs R_T^{1/2}\log(e+\|\g u\|_2)  + R_T\\
	&\leqs R_T^{1/2}\Big(C(\a)R_T^{(1+\a\beta)/2} + R_T^{(\gamma-2\beta+1)/2}\Big)  +1\\
	&\leqs C(\a)R_T^{1+\a\beta} + R_T^{1+\gamma/2-\beta}\\
	\ea\]

\endproof

Next, we are ready to give the upper bound of the density independent of the time T.
\begin{Proposition}\label{SZ1} Under the conditions of Theorem \ref{them1}, and 
\be
\max\{1,\frac{\gamma+2}{4}\}<\beta\le\gamma.
\ee
then there is a constant $C>0$ independent of the time T, such that
	~\be
	\parallel \rho \parallel_{L^{\infty} L^{\infty}}  \leq C.
	\ee
\end{Proposition}
\proof
Recall \eqref{thxi} that
\be\label{thxi-1}
(\theta +\xi )_{t} + u \cdot \partial_{r}(\theta+\xi) +\int_{R}^{r} \frac{\rho u^2}{s} ds +P + F(R,\ t)=0.
\ee

The main idea is to employ Zlotnik's inequality to get the time-independent estimate for the density. We just need to rewrite the equation in suitable form and verify the conditions required by the Zloinik type inequality. This will result in the restriction of $\beta$.

Set 
\be
y=\theta+\xi,\ g(y)=-P(\theta^{-1}(y-\xi)),\ g(+\infty)=-\infty,
\ee
and
\be
h=\int_0^t\int_R^r \f{\r u^2}{s}dsd\tau + \int_0^t F(R,\tau)d\tau.
\ee
Then \eqref{thxi-1} can be written as
\be
\frac{D}{Dt}y=g(y)+h_t.
\ee
It suffices to verifies the conditions required in lemma \ref{zlemma}, especially the Lipschitz coefficients of $h$.

For any time $0<t_1<t_2$,
\[
h(t_2)-h(t_1)=\int_{t_1}^{t_2}\int_R^r \f{\r u^2}{s}dsd\tau + \int_{t_1}^{t_2} F(R,\tau)d\tau
\]

\[\ba
|h(t_2)-h(t_1)|&=\int\int_R^r \f{\r u^2}{s}dsd\tau + \int_{t_1}^{t_2} F(R,\tau)d\tau\\
&\leqs \|\rho\|_\infty\int_0^T\|\g u\|_2^2dt +\Big|\int_{t_1}^{t_2} F(R,\tau)d\tau\Big|\\
&\leq C\|\rho\|_\infty + C + C(t_2-t_1).
\ea\]
where we had used the following fact
\be\label{fr-1-2}
\ba
\Big|\int_{t_1}^{t_2} F(R,\tau)d\tau\Big| \leq C\|\rho\|_\infty + C + C(t_2-t_1).
\ea
\ee
 
To prove \eqref{fr-1-2}, we make use of Lemma \ref{FR}. Indeed, we had
\be\label{fr-1-1}
\ba
\int_{t_1}^{t_2} F(R,\tau)d\tau &=  \frac{1}{R^2} \Big[\int_{0}^{R} \rho u r^2 dr \big|_{\ts=t_2} + \int_{0}^{R} \rho u r^2 dr \big|_{\ts=t_1}\\
& + 2\int_{t_1}^{t_2}\int_{0}^{R} F(r, \ts) r dr d\ts - \int_{t_1}^{t_2}\int_{0}^{R} \rho u^2 r dr d\ts \Big]
\ea
\ee
Consequently,

\[\ba
\Big|\int_{t_1}^{t_2} F(R,\tau)d\tau\Big| &\leqs 1 + \Big|\int_{t_1}^{t_2}\int_{0}^{R} F(r, \ts) r dr d\ts\Big| +\| \rho\|_\infty\int_0^T\|\div u\|_2^2dt\\
&\leqs 1 + \| \rho\|_\infty + \Big|\int_{t_1}^{t_2}\int_{0}^{R} F(r, \ts) r dr d\ts\Big|.
\ea\]
Note that
\[\ba
\Big|\int_{0}^{R} F(r, \ts) r dr\Big| &= \Big|\int F(x, \ts) dx \Big|=\Big|\int \ld\div u - p dx\Big|\\
&\lesssim \int \ld + \ld|\div u|^2 + p dx\Big|\\
&\lesssim \int \ld dx + \int \ld|\div u|^2dx + 1,
\ea\]
which gives that
\[
\Big|\int_{t_1}^{t_2}\int_{0}^{R} F(r, \ts) r dr d\ts \Big|
\leqs \int_{t_1}^{t_2}\int \ld dx + \int \ld|\div u|^2dx + 1d\ts
\leq C+C(t_2-t_1).
\]
with the help of the energy inequality and $\beta\le\gamma$.

Hence we arrive at
\[\ba
\Big|\int_{t_1}^{t_2} F(R,\tau)d\tau\Big| \leq C\|\rho\|_\infty + C + C(t_2-t_1).
\ea\]

By Zlotnik's lemma,
\[
\theta + \xi \leqs \|\rho\|_\infty +1.
\]
\[
\rho^\beta \leqs \theta \leqs \|\xi\|_\infty + \|\rho\|_\infty +1 \leqs C(\a)R_T^{1+\a\beta} + R_T^{1+\gamma/2-\beta}.
\]
then,
\[
R_T^{\beta} \leqs C(\a)R_T^{1+\a\beta} + R_T^{1+\gamma/2-\beta}.
\]
Direct calculations show that
\be
R_T\le C
\ee
independent of time $T$ as long as 
\be
\max\{1,\frac{\gamma+2}{4}\}<\beta\le\gamma.
\ee
This finishes the proof of Proposition \ref{SZ1}.
\endproof

\section{Proofs of Theorem \ref{them1} and Theorem \ref{them2}}
With all the a priori estimates at hand, the proofs of Theorem \ref{them1} and Theorem \ref{them1} are standard, so omited.
See \cite{L2021a,2016Huang} for details.

\section*{Acknowledgements}
 X.-D. Huang is partially supported by NNSFC Grant Nos. 11971464, 11688101 and CAS Project for Young Scientists in Basic Research, Grant No.YSBR-031, National Key R\&D Program of China, Grant No.2021YFA1000800. W. YAN is partially supported by NNSFC Grant Nos. 11371064, 11871113. R. Yu is partially supported by Guangdong Provincial Natural Science Foundation, Grant No. 2020A1515110942.

\end{document}